\documentclass[12pt]{amsart}
\usepackage{amsxtra,latexsym,amssymb}
\usepackage{epsfig,rotate,amsthm}
\usepackage{hyperref}

\def\qed{{\hfill $\Box$}}

\def\Z{{\mathbb Z}}
\def\C{{\mathbb C}}

\def\U{{U_{r,s}^{+}(G_{2})}}
\def\V{{\check U}^{\geq 0}_{r,s}(G_{2})}
\theoremstyle{theorem}
\newtheorem{thm}{Theorem}[section]
\newtheorem{cor}{Corollary}[section]
\newtheorem{prop}{Proposition}[section]

\newtheorem{lem}{Lemma}[section]
\theoremstyle{definition}
\newtheorem{defn}{Definition}[section]
\theoremstyle{remark}

\begin{document}
\title[Automorphisms of $\V$]{Automorphisms of the two-parameter Hopf algebra $\V$}
\author[X. Tang]{Xin Tang}
\address{Department of Mathematics \& Computer Science\\
Fayetteville State University\\
1200 Murchison Road, Fayetteville, NC 28301}
\email{xtang@uncfsu.edu} 
\keywords{Algebra automorphisms, Hopf algebra automorphisms, Two-parameter quantized enveloping algebra}
\thanks{}
\date{\today}
\subjclass[2000]{Primary 17B37,16B30,16B35.}
\begin{abstract}
We determine the group of algebra automorphisms for the two-parameter quantized enveloping algebra $\V$. As an application, we prove that the group of Hopf algebra automorphisms for 
$\V$ is isomorphic to a torus of rank two. 
\end{abstract}
\maketitle
\section*{Introduction}

Let $\mathfrak g$ be a finite dimensional complex simple Lie algebra and let $r, s\in \C^{\ast}$. The two-parameter quantized enveloping algebras (or quantum groups) $U_{r,s}({\mathfrak g})$ have been studied in the literatures \cite{BW1, BW2, BGH} and the references therein. Recently, more studies have been conducted toward their subalgebras such as $U_{r,s}^{+}(\mathfrak g)$,  and the augmented Hopf algebras $\check{U}_{r,s}^{\geq 0}(\mathfrak g)$. In \cite{T1}, the author has computed the derivations for the subalgebra $U_{r,s}^{+}(sl_{3})$, and determined both the algebra automorphism group and Hopf algebra automorphism group for the Hopf algebra $\check{U}_{r,s}^{\geq 0}(sl_{3})$.  A similar work has been carried out for the algebra $U_{r,s}^{+}(B_{2})$ and the Hopf algebra $\check{U}_{r,s}^{\geq 0}(B_{2})$ in \cite{T2}. The results in these works suggest that the subalgebras $U_{r,s}^{+}(\mathfrak g)$ and $\check{U}_{r,s}^{\geq 0}(\mathfrak g)$ are close analogues of 
their one-parameter counterparts, which facilitates further investigation toward these subalgebras.

In this paper, we are planning to derive some similar results for the two-parameter Hopf algebra $\V$ in terms of its (Hopf) algebra automorphisms. In particular, we will first determine the group of 
algebra automorphisms for the Hopf algebra  $\V$. Then, as an application, we will further prove that the group of Hopf algebra automorphisms  for$\V$ is indeed isomorphic to a torus of rank 2. We 
will closely follow the approach used in \cite{F}.

The paper is organized as follows. In Section 1, we will recall some basics on the two-parameter Hopf algebra $\V$. In Section 2, we will determine the group of algebra automorphisms for the Hopf algebra $\V$. In Section 3, we will determine the group of Hopf algebra automorphisms for the Hopf algebra $\V$. 

\section{Some basic properties of the Hopf algebra $\V$}

Recall that the two-parameter quantum group $U_{r,s}(G_{2})$ associated to the finite dimensional complex simple Lie algebra of type $G_{2}$ has been studied in \cite{HQ,HW}. In particular, a PBW basis of $U_{r,s}(G_{2})$ has been constructed in \cite{HW}. For the readers{'} convenience, we will recall the construction of the subalgebra $\U$ together with some of its basic properties from \cite{HQ}. In the rest of this paper, we will always assume that the parameters $r, s$ are chosen from $\C^{\ast}$ such that $r^{m}s^{n}=1$ implies $m=n=0$. 

First of all, we need to recall the following definition from the references \cite{HQ, HW}:
\begin{defn}
The two-parameter quantized enveloping algebra $\U$ is defined to be the $\C-$algebra 
generated by the generators $e_{1}, e_{2}$ subject to the following relations:
\begin{eqnarray*}
e_{2}^{2}e_{1}-(r^{-3}+s^{-3})e_{2}e_{1}e_{2}+r^{-3}s^{-3}e_{1}e_{2}^{2}=0,\\
e_{1}^{4}e_{2}-(r+s)(r^{2}+s^{2})e_{1}^{3}e_{2}e_{1}+rs(r^{2}+rs+s^{2})e_{1}^{2}e_{2}e_{1}^{2}\\-r^{3}s^{3}(r+s)(r^{2}+s^{2})e_{1}e_{2}e_{2}^{3}+r^{6}s^{6}e_{2}e_{1}^{4}=0.
\end{eqnarray*}
\end{defn}

In the rest of this section, we will establish some basic properties of the algebra $\U$ and introduce an augmented Hopf algebra $\V$. In particular, we will recall the construction of a PBW basis for the algebra $U$.

In order to recall the construction of a PBW basis of $\U$, we  need to fix some variables. We will follow the notation in \cite{HW}.
\begin{eqnarray*}
X_{6}=E_{1}=e_{1},\quad X_{1}=E_{2}=e_{2},\\
X_{2}=E_{12}=e_{1}e_{2}-s^{3}e_{2}e_{1},\\
X_{4}=E_{112}=e_{1}E_{12}-rs^{2}E_{12}e_{1},\\
X_{5}=E_{1112}=e_{1}E_{112}-r^{2}sE_{112}e_{1},\\
X_{3}=E_{11212}=E_{112}E_{12}-r^{2}sE_{12}E_{112}.
\end{eqnarray*}

Now  we can have the following result
\begin{thm}
({\bf Theorem 2.4.} in \cite{HW}) The following set 
\[
\{ E_{2}^{n_{1}}E_{12}^{n_{2}}E_{11212}^{n_{3}}
E_{112}^{n_{4}}E_{1112}^{n_{5}}E_{1}^{n_{6}}\mid n_{i}\in \Z_{\geq 0} \}
\]
forms a Lyndon basis of  the algebra $\U$.
\end{thm}
\qed

We now recall the definition of the Hopf subalgebra $U^{\geq 0}_{r,s}(G_{2})$ from \cite{HQ, HW}. We shall have the following definition. 
\begin{defn}
The Hopf algebra $U_{r,s}^{\geq 0}(G_{2})$ is defined to be the $\C-$algebra generated by the generators $e_{1}, e_{2}$ and $w_{1}, w_{2}$ subject to the following relations:
\begin{eqnarray*}
w_{1}w_{1}^{-1}=w_{2}w_{2}^{-1}=1,\quad w_{1}w_{2}=w_{2}w_{1};\\
w_{1}e_{1}=rs^{-1}e_{1}w_{1},\quad w_{1}e_{2}=s^{3}e_{2}w_{1};\\
w_{2}e_{1}=r^{-3}e_{1}w_{2},\quad w_{2}e_{2}=r^{3}s^{-3}e_{2}w_{2};\\
e_{2}^{2}e_{1}-(r^{-3}+s^{-3})e_{2}e_{1}e_{2}+r^{-3}s^{-3}e_{1}e_{2}^{2}=0;\\
e_{1}^{4}e_{2}-(r+s)(r^{2}+s^{2})e_{1}^{3}e_{2}e_{1}+rs(r^{2}+rs+s^{2})e_{1}^{2}e_{2}e_{1}^{2}\\-r^{3}s^{3}(r+s)(r^{2}+s^{2})e_{1}e_{2}e_{2}^{3}+r^{6}s^{6}e_{2}e_{1}^{4}=0.
\end{eqnarray*}
\end{defn}

Let us set
\begin{eqnarray*}
\Delta(e_{1})=e_{1}\otimes 1+w_{1}\otimes e_{1};\\
\Delta(e_{2})=e_{2}\otimes 1+w_{2}\otimes e_{2};\\
\Delta(w_{1})=w_{1}\otimes w_{1},\quad \Delta(w_{2})=w_{2}\otimes w_{2};\\
S(e_{1})=-w_{1}e_{1},\quad S(e_{2})=-w_{2}e_{2};\\
S(w_{1})=w_{1}^{-1},\quad S(w_{2})=w_{2}^{-1};\\
\epsilon(e_{1})=\epsilon(e_{2})=0,\quad \epsilon (w_{1})=\epsilon(w_{2})=1.
\end{eqnarray*}
Then, it is easy to see that the above operators define a Hopf algebra structure on the $U^{\geq 0}_{r,s}(G_{2})$; and we further have the following proposition:
\begin{prop}
The set 
\[
\{X_{1}^{a}X_{2}^{b}X_{3}^{c}X_{4}^{d}X_{5}^{e}X_{6}^{f}w_{1}^{m}w_{2}^{n}|a, b, c, d , e, f \in \Z_{\geq 0}, m, n \in \Z\}
\]
forms a PBW-basis of the Hopf algebra $U_{r,s}^{\geq 0}(G_{2})$.
\end{prop}
\qed

To define the augmented Hopf algebra $\V$, let us set the following new variables
\[
k_{1}=w_{1}^{2}w_{2},\quad k_{2}=w_{1}^{3}w_{2}^{2}.
\]

Now we have the following definition of the augmented Hopf algebra $\V$.
 \begin{defn}
The Hopf algebra $\V$ is defined to be the $\C-$algebra generated by the generators $e_{1}, e_{2}$ and $k_{1}, k_{2}$ subject to the following relations:
\begin{eqnarray*}
k_{1}k_{1}^{-1}=k_{2}k_{2}^{-1}=1,\quad k_{1}k_{2}=k_{2}k_{1};\\
k_{1}e_{1}=r^{-1}s^{-2}e_{1}k_{1},\quad k_{1}e_{2}=r^{3}s^{3}e_{2}k_{1};\\
k_{2}e_{1}=r^{-3}s^{-3}e_{1}k_{2},\quad k_{2}e_{2}=r^{6}s^{3}e_{2}k_{2};\\
e_{2}^{2}e_{1}-(r^{-3}+s^{-3})e_{2}e_{1}e_{2}+r^{-32}s^{-3}e_{1}e_{2}^{2}=0;\\
e_{1}^{4}e_{2}-(r+s)(r^{2}+s^{2})e_{1}^{3}e_{2}e_{1}+rs(r^{2}+rs+s^{2})e_{1}^{2}e_{2}e_{1}^{2}\\-r^{3}s^{3}(r+s)(r^{2}+s^{2})e_{1}e_{2}e_{2}^{3}+r^{6}s^{6}e_{2}e_{1}^{4}=0.
\end{eqnarray*}
\end{defn}

Once again, let us further define the following
\begin{eqnarray*}
\Delta(e_{1})=e_{1}\otimes 1+k_{1}^{2}k_{2}^{-1}\otimes e_{1};\\
\Delta(e_{2})=e_{2}\otimes 1+k_{1}^{-3}k_{2}^{2}\otimes e_{2};\\
\Delta(k_{1})=k_{1}\otimes k_{1},\quad \Delta(k_{2})=k_{2}\otimes k_{2};\\
S(e_{1})=-k_{1}^{2}k_{2}^{-1}e_{1},\quad S(e_{2})=-k_{1}^{-3}k_{2}^{2}e_{2};\\
S(k_{1})=k_{1}^{-1},\quad S(k_{2})=k_{2}^{-1};\\
\epsilon(e_{1})=\epsilon(e_{2})=0,\quad \epsilon (k_{1})=\epsilon(k_{2})=1.
\end{eqnarray*}

Then, it is easy to see that the above operators define a Hopf algebra structure on the augmented Hopf algebra $\V$. Furthermore, the Hopf algebra $\V$ has a PBW basis as described below.
\begin{prop}
The set 
\[
\{X_{1}^{a}X_{2}^{b}X_{3}^{c}X_{4}^{d}X_{5}^{e}X_{6}^{f}k_{1}^{m}k_{2}^{n}|a, b, c, d , e, f \in \Z_{\geq 0}, m, n \in \Z\}
\]
forms a basis of $\V$ over the base field $\C$.
\end{prop}  
\qed

\section{Algebra automorphism group of the Hopf algebra $\V$}
In this section, we will determine the algebra automorphism group of the Hopf algebra $\V$. We 
will closely follow the approach used in \cite{F}. Note that such an approach has been adopted to investigate the automorphism group of the Hopf algebra $\check{U}_{r,s}^{\geq 0}(sl_{3})$ in \cite{T1}. Similar work has also appeared in \cite{T2}. It is no surprise that we will derive very similar results to those obtained in \cite{T1, T2}.

Let us denote by $\theta$ an algebra automorphism of the Hopf algebra $\V$. Since the elements $k_{1}, k_{2}$ are invertible elements in the algebra $\V$ and $\theta$ is an algebra automorphism of the Hopf algebra  $\V$, the images $\theta(k_{1}), \theta(k_{2})$  of the invertible elements $k_{1}, k_{2}$ are invertible elements in the Hopf algebra $\V$. Note that it 
is easy to check  that the only invertible elements of the algebra $\V$ are of the form $\lambda k_{1}^{m}k_{2}^{n}, \lambda \in \C^{\ast}, m, n \in \Z$. Therefore, the elements $\theta(k_{1})$ and $\theta(k_{2})$ can be expressed as follows
\[
\theta(k_{1})=\lambda_{1}k_{1}^{x}k_{2}^{y},\quad \theta(k_{2})=\lambda_{2}k_{1}^{z}k_{2}^{w}
\]
for some $\lambda_{1}, \lambda_{2}\in \C^{\ast}$ and some $x,y,z,w \in \Z$. 

Note that we can also associate an invertible $2\times 2$ matrix to the algebra automorphism $\theta$; and we will denote this matrix by $M_{\theta}=(M_{ij})$. As a matter of fact, we 
will define this matrix by the entries as follows
\[
M_{11}=x, \, M_{12}=y, \, M_{21}=z, \, M_{22}=w.
\]

Due to the fact that the mapping $\theta$ is an algebra automorphism, the determinant of 
the corresponding matrix $M_{\theta}$ is either $1$ or $-1$. That is, we shall have that 
\[
xw-yz=\pm 1.
\]

In terms of the PBW basis of $\V$, we can further express the images of the generators $e_{1}, e_{2}$ of $\V$ under the algebra automorphism $\theta$ as follows
\[
\theta(e_{l})=\sum_{m_{l}, n_{l}, \beta_{l}^{1}, \beta_{l}^{2}, \beta_{l}^{3},\beta_{l}^{4}, \beta_{l}^{5},\beta_{l}^{6}}\gamma_{m_{l} n_{l}\beta_{l}^{1}\beta_{l}^{2}\beta_{l}^{3}\beta_{l}^{4}\beta_{l}^{5}\beta_{l}^{6}}k_{1}^{m_{l}}k_{2}^{n_{l}}X_{1}^{\beta_{l}^{1}}X_{2}^{\beta_{l}^{2}}X_{3}^{\beta_{l}^{3}} X_{4}^{\beta_{l}^{4}}X_{5}^{\beta_{l}^{5}}\beta_{6}^{\beta_{l}^{6}}
\]
where $\gamma_{m_{l}n_{l}\beta_{l}^{1}\beta_{l}^{2}\beta_{l}^{3} \beta_{l}^{4}}\beta_{l}^{5}\beta_{l}^{6}$ are chosen from $\C^{\ast}$, and $m_{l}, n_{l}$ are chosen from  $\Z$, and $\beta_{l}^{1}, \beta_{l}^{2}, \beta_{l}^{3}, \beta_{l}^{4}, \beta_{l}^{5}$ and $\beta_{l}^{6}$ are chosen from  $\Z_{\geq 0}$. 

In the rest of this section, we prove that $\theta$ is actually defined in a simple and specific way. First of all, we are going to establish some identities, whose proofs involve straightforward verifications; and we will not record these verifications.
 
\begin{lem} For $l=1, 2$, the following identities shall hold
\begin{eqnarray*}
&&k_{1}^{x}k_{2}^{y}X_{1}^{\beta_{l}^{1}}X_{2}^{\beta_{l}^{2}}
X_{3}^{\beta_{l}^{3}}X_{4}^{\beta_{l}^{4}}X_{5}^{\beta_{l}^{5}}X_{6}^{\beta_{l}^{6}}\\
&=&
(r^{-1})^{x(-3\beta_{l}^{1}-2\beta_{l}^{2}-3\beta_{l}^{3}-\beta_{l}^{4}+\beta_{l}^{6})
+y(-6\beta_{l}^{1}-3\beta_{l}^{2}-3\beta_{l}^{3}+3\beta_{l}^{5}+3\beta_{l}^{6})}\\
&&(s^{-2})^{x(-3\beta_{l}^{1}-2\beta_{l}^{2}-3\beta_{l}^{3}-\beta_{l}^{4}+\beta_{l}^{6})
+y(-6\beta_{l}^{1}-3\beta_{l}^{2}-3\beta_{l}^{3}+3\beta_{l}^{5}+3\beta_{l}^{6})}\\
&&X_{1}^{\beta_{l}^{1}}X_{2}^{\beta_{l}^{2}}
X_{3}^{\beta_{l}^{3}}X_{4}^{\beta_{l}^{4}}X_{5}^{\beta_{l}^{5}}
X_{6}^{\beta_{l}^{6}}k_{1}^{x}k_{2}^{y}.
\end{eqnarray*}
\end{lem}
\qed

Now we have the following proposition, which characterizes the nature of the matrix $M_{\theta}$.
\begin{prop}
Let $\theta \in Aut_{\C}(\V)$ be an algebra automorphism of the Hopf algebra $\V$, then we 
have $M_{\theta} \in GL(2,\Z_{\geq 0})$.
\end{prop}
{\bf Proof:} Since $k_{1}e_{1}=r^{-1}s^{-2}e_{1}k_{1}, k_{2}e_{1}=r^{-1}s^{-1}e_{1}k_{2}$ and $\theta$ is an algebra automorphism of $\V$, we have the following
\begin{eqnarray*}
\theta(k_{1})\theta(e_{1})=r^{-1}s^{-2}\theta(e_{1})\theta(k_{1});\\
\theta(k_{2})\theta(e_{1})=r^{-3}s^{-3}\theta(e_{1})\theta(k_{2}).
\end{eqnarray*}

Using the previous lemma, we shall have the following identities
\begin{eqnarray*}
 x(-3\beta_{1}^{1}-2\beta_{1}^{2}-3\beta_{1}^{3}-\beta_{1}^{4}+
\beta_{1}^{6})+y(-6\beta_{11}^{1}-3\beta_{1}^{2}-3\beta_{1}^{3}+
3\beta_{1}^{5}+3\beta_{1}^{6})\\
=1;\\
x(-3\beta_{1}^{1}-2\beta_{1}^{2}-3\beta_{1}^{3}-\beta_{1}^{4}+\beta_{1}^{6})+
y(-6\beta_{11}^{1}-3\beta_{1}^{2}-3\beta_{1}^{3}+3\beta_{1}^{5}+3\beta_{1}^{6})
\\=2;\\
x(-3\beta_{1}^{1}-2\beta_{1}^{2}-3\beta_{1}^{3}-\beta_{1}^{4}+\beta_{1}^{6})+
y(-6\beta_{11}^{1}-3\beta_{1}^{2}-3\beta_{1}^{3}+3\beta_{1}^{5}+3\beta_{1}^{6})\\
=-3; \\
x(-3\beta_{1}^{1}-2\beta_{1}^{2}-3\beta_{1}^{3}-\beta_{1}^{4}+\beta_{1}^{6})+
y(-6\beta_{11}^{1}-3\beta_{1}^{2}-3\beta_{1}^{3}+3\beta_{1}^{5}+3\beta_{1}^{6})\\
=-3;\\
z(-3\beta_{1}^{1}-2\beta_{1}^{2}-3\beta_{1}^{3}-\beta_{1}^{4}+\beta_{1}^{6})+
w(-6\beta_{11}^{1}-3\beta_{1}^{2}-3\beta_{1}^{3}+3\beta_{1}^{5}+3\beta_{1}^{6})\\
=3;\\
z(-3\beta_{1}^{1}-2\beta_{1}^{2}-3\beta_{1}^{3}-\beta_{1}^{4}+\beta_{1}^{6})+
w(-6\beta_{11}^{1}-3\beta_{1}^{2}-3\beta_{1}^{3}+3\beta_{1}^{5}+3\beta_{1}^{6})\\
=3;\\z(-3\beta_{1}^{1}-2\beta_{1}^{2}-3\beta_{1}^{3}-\beta_{1}^{4}+\beta_{1}^{6})+
w(-6\beta_{11}^{1}-3\beta_{1}^{2}-3\beta_{1}^{3}+3\beta_{1}^{5}+3\beta_{1}^{6})\\
=-6;\\z(-3\beta_{1}^{1}-2\beta_{1}^{2}-3\beta_{1}^{3}-\beta_{1}^{4}+\beta_{1}^{6})+
w(-6\beta_{11}^{1}-3\beta_{1}^{2}-3\beta_{1}^{3}+3\beta_{1}^{5}+3\beta_{1}^{6})\\
=-3.
\end{eqnarray*}

After some combinations and simplifications of these equations,  we shall have the following system of equations:
\begin{eqnarray*}
x(\beta_{1}^{2}+3\beta_{1}^{3}+2\beta_{1}^{4}+3\beta_{1}^{5}+\beta_{1}^{6})+
y(3\beta_{1}^{1}+3\beta_{1}^{2}+6\beta_{1}^{3}+3\beta_{1}^{4}+3\beta_{1}^{5})\\
=1;\\
x(\beta_{2}^{2}+3\beta_{2}^{3}+2\beta_{2}^{4}+3\beta_{2}^{5}+\beta_{2}^{6})+
y(3\beta_{2}^{1}+3\beta_{2}^{2}+6\beta_{2}^{3}+3\beta_{2}^{4}+3\beta_{2}^{5})\\
=0;\\
z(\beta_{1}^{2}+3\beta_{1}^{3}+2\beta_{1}^{4}+3\beta_{1}^{5}+\beta_{1}^{6})+
w(3\beta_{1}^{1}+3\beta_{1}^{2}+6\beta_{1}^{3}+3\beta_{1}^{4}+3\beta_{1}^{5})\\
=0;\\
z(\beta_{2}^{2}+3\beta_{2}^{3}+2\beta_{2}^{4}+3\beta_{2}^{5}+\beta_{2}^{6})+
w(3\beta_{2}^{1}+3\beta_{2}^{2}+6\beta_{2}^{3}+3\beta_{2}^{4}+3\beta_{2}^{5})\\
=3.
\end{eqnarray*}

Now let us define a $2\times2-$matrix $B=(b_{ij})$ by setting the entries of $B$ as follows: 
\begin{eqnarray*}
b_{11}=\beta_{1}^{2}+3\beta_{1}^{3}+2\beta_{1}^{4}+3\beta_{1}^{5}+\beta_{1}^{6};\\
b_{21}=3\beta_{1}^{1}+3\beta_{1}^{2}+6\beta_{1}^{3}+3\beta_{1}^{4}+3\beta_{1}^{5};\\
b_{12}=\beta_{1}^{2}+3\beta_{1}^{3}+2\beta_{1}^{4}+3\beta_{1}^{5}+\beta_{1}^{6};\\
b_{22}=3\beta_{2}^{1}+3\beta_{2}^{2}+6\beta_{2}^{3}+3\beta_{2}^{4}+3\beta_{2}^{5}.
\end{eqnarray*}

Thus we shall have the following
\[
M_{\theta}B=\left( \begin{array}{lr}
1 & 0\\
0 &3 
\end{array}
\right)
\]
which implies that we have  
\[
M_{\theta}^{-1}=\left( \begin{array}{lr}
b_{11} & b_{12}/3\\
b_{21} &b_{22} /3
\end{array}
\right).
\] 

Let us denote by $M_{\theta^{-1}}$ the matrix associated to the inverse of $\theta$, then we have $M_{\theta^{-1}}=M_{\theta}^{-1}$. Since all the entries $b_{11}, b_{12}, b_{21}, b_{22}$ of the matrix $B$ are all nonnegative integers, we know that the matrix $M_{\theta^{-1}}$ is indeed in the group $GL(2, \Z_{\geq 0})$. Apply this process to the 
algebra automorphism $\theta^{-1}$, we have that the matrix $M_{\theta}$ is in $GL(2, \Z_{\geq 0})$. So we have proved the proposition.
\qed

In addition, please note that the following  important lemma was already established in the reference \cite{F}. This lemma applies to our case as well. 
\begin{lem}
If $M$ is a matrix in $GL(n,\Z_{\geq 0})$ such that its inverse matrix $M^{-1}$ is also in $GL(n,\Z_{\geq 0})$, then we 
have $M=(\delta_{i\sigma(j)})_{i,j}$, where $\sigma$ is an element of the symmetric group $\mathbb{S}_{n}$.
\end{lem}
\qed

Based on {\bf Proposition 2.1} and {\bf Lemma 2.2}, it is easy to see that we have the following result.
\begin{cor}
Suppose that $\theta \in Aut_{\C}(\V)$ is an algebra automorphism of the Hopf algebra 
$\V$. Then for $l=1,2$, we have 
\[
\theta(k_{l})=\lambda_{l} k_{\sigma(l)}
\]
where $\sigma \in \mathbb{S}_{2}$ and $\lambda_{l} \in \C^{\ast}$.
\end{cor}
\qed

To proceed, we need some further preparations. Suppose that we have $\theta(k_{1})=\lambda_{1}k_{1}$ and $\theta(k_{2})=\lambda_{2}k_{2}$. Then 
we have the following 
\begin{lem}
The following identities hold
\begin{eqnarray*}
-3\beta_{1}^{1}-2\beta_{1}^{2}-3\beta_{1}^{3}-\beta_{1}^{4}+\beta_{1}^{6}=1;\\
-3\beta_{1}^{1}-\beta_{1}^{2}+\beta_{1}^{4}+3\beta_{1}^{5}+\beta_{1}^{6}=2;\\
-6\beta_{1}^{1}-3\beta_{1}^{2}-3\beta_{1}^{3}+3\beta_{1}^{5}+3\beta_{1}^{6}=3;\\
-3\beta_{1}^{1}+3\beta_{1}^{3}+3\beta_{1}^{4}+6\beta_{1}^{5}+3\beta_{1}^{6}=3;\\
-3\beta_{2}^{1}-2\beta_{2}^{2}-3\beta_{2}^{3}-\beta_{2}^{4}+\beta_{2}^{6}=-3;\\
-3\beta_{2}^{1}-\beta_{2}^{2}+\beta_{2}^{4}+3\beta_{2}^{5}+\beta_{2}^{6}=-3;\\
-6\beta_{2}^{1}-3\beta_{2}^{2}-3\beta_{2}^{3}+3\beta_{2}^{5}+3\beta_{2}^{6}=-6;\\
-3\beta_{2}^{1}+3\beta_{2}^{3}+3\beta_{2}^{4}+6\beta_{2}^{5}+3\beta_{2}^{6}=-3.
\end{eqnarray*}
\end{lem}
\qed

Moreover, the identities in the previous lemma imply the following
\begin{lem}
The following identities hold
\begin{eqnarray*}
\beta_{1}^{2}+3\beta_{1}^{3}+2\beta_{1}^{4}+3\beta_{1}^{5}+\beta_{1}^{6}=1;\\
3\beta_{1}^{1}+3\beta_{1}^{2}+6\beta_{1}^{3}+3\beta_{1}^{5}=0;\\
\beta_{2}^{2}+3\beta_{2}^{3}+2\beta_{2}^{4}+3\beta_{2}^{5}+\beta_{2}^{6}=1;\\
3\beta_{2}^{1}+3\beta_{2}^{2}+6\beta_{2}^{3}+3\beta_{2}^{5}=0.\\
\end{eqnarray*}
In particular, we have the following
\begin{eqnarray*}
\beta_{1}^{1}=\beta_{1}^{2}=\beta_{1}^{3}=\beta_{1}^{4}=\beta_{1}^{5}=0,\quad \beta_{1}^{6}=1;\\
\beta_{2}^{2}=\beta_{2}^{3}=\beta_{2}^{4}=\beta^{5}=\beta_{2}^{6}=0,\quad \beta_{2}^{1}=1.
\end{eqnarray*}
\end{lem}
\qed

Similarly, if we assume that we have $\theta(k_{1})=\lambda_{1}k_{2}$ and $\theta(k_{2})=\lambda_{2}k_{1}$, then we shall have the following
\begin{lem}
\begin{eqnarray*}
 \beta_{1}^{2}=\beta_{1}^{3}=\beta_{1}^{4}=\beta_{1}^{5}=\beta_{1}^{6}=0,\quad \beta_{1}^{1}=1;\\
\beta_{2}^{1}=\beta_{2}^{2}=\beta_{2}^{3}=\beta^{4}=\beta_{2}^{5}=0,\quad \beta_{2}^{6}=1.
\end{eqnarray*}
\end{lem}
\qed

Follows from the previous two lemmas, we can easily have the following result
\begin{prop}
Let $\theta \in Aut_{\C}(\V)$ be an algebra automorphism of the Hopf algebra $\V$. Then for $l=1,2$, we have 
\[
\theta(e_{l})=\gamma_{l}k_{1}^{m_{l}}k_{2}^{n_{l}}e_{\sigma(l)}
\]
where $\gamma_{l}\in \C^{\ast}$ and $m_{l}, n_{l} \in \Z$.
\end{prop}
\qed

The following result will further demonstrate that the two generators $e_{1}, e_{2}$ of $\V$ can not be exchanged by any algebra automorphism $\theta$ of $\V$. In particular, we have the following result
\begin{cor}
Let $\theta \in Aut_{\C}(\V)$ be an algebra automorphism of $\V$. Then for $l=1,2$, we have the following 
\[
\theta(k_{l})=\lambda_{l}k_{l},\, \theta(e_{l})=\gamma_{l}k_{1}^{m_{l}}k_{2}^{n_{l}}e_{l}
\]
where $\lambda_{l}, \gamma_{l} \in \C^{\ast}$ and $m_{l},n_{l} \in \Z$.
\end{cor}
{\bf Proof:} Suppose that $\theta(k_{1})=\lambda_{1}k_{2}$ and $\theta(e_{2})=\gamma_{1}k_{1}^{m_{1}}k_{2}^{n_{1}}e_{2}$. Since we have $\theta(k_{1})\theta(e_{1})=r^{-1}s^{-2}\theta(e_{1})\theta(k_{1})$, we have the 
following
\[
\lambda_{1}k_{2}\gamma_{1} k_{1}^{m_{1}}k_{2}^{m_{2}}e_{2}=r^{-1}s^{-2}\gamma_{1} k_{1}^{m_{1}}k_{2}^{m_{2}}e_{2}\lambda_{1}k_{2}.
\]
Note that $k_{2}e_{2}=r^{6}s^{3}e_{2}k_{2}$, then we got a contradiction. Therefore, we have proved the statement as desired.
\qed

Now we will further establish some identities via direct calculations and we will skip the detailed calculations here.

\begin{lem} We have the following identities:
\begin{eqnarray*}
(k_{1}^{a}k_{2}^{b}e_{1})^{4}(k_{1}^{c}k_{2}^{d}e_{2})&=&r^{6a+18b+4c+12d}s^{12a+18b+8c+12d}
k_{1}^{a+c}k_{2}^{b+d}e_{1}^{4}e_{2};\\
(k_{1}^{a}k_{2}^{b}e_{1})^{3}(k_{1}^{c}k_{2}^{d}e_{2})(k_{1}^{a}k_{2}^{b}e_{1})&=&r^{3a+12b+3c+9d}s^{9a+15b+6c+9d}
k_{1}^{4a+c}k_{2}^{4b+d}e_{1}^{3}e_{2}e_{1};\\
(k_{1}^{a}k_{2}^{b}e_{1})^{2}(k_{1}^{c}k_{2}^{d}e_{2})(k_{1}^{a}k_{2}^{b}e_{1})^{2})&=&r^{6b+2c+6d}s^{6a+12b+4c+6d}
k_{1}^{4a+c}k_{2}^{4b+d}e_{1}^{2}e_{2}e_{1}^{2};\\
(k_{1}^{a}k_{2}^{b}e_{1})(k_{1}^{c}k_{2}^{d}e_{2})(k_{1}^{a}k_{2}^{b}e_{1})^{3}&=&r^{-3a+c+3d}s^{3a+9b+2c+3d}
k_{1}^{4a+c}k_{2}^{4b+d}e_{1}e_{2}e_{1}^{3};\\
(k_{1}^{c}k_{2}^{d}e_{2})(k_{1}^{a}k_{2}^{b})^{4}&=&r^{-6a-6b}s^{6b}
k_{1}^{4a+c}k_{2}^{4b+d}e_{2}e_{1}^{4}.
\end{eqnarray*}
\end{lem}
\qed

Similarly, we can also have the following lemma, whose proof will be skipped.
\begin{lem} The following identities hold.
\begin{eqnarray*}
(k_{1}^{c}k_{2}^{d}e_{2})^{2}(k_{1}^{a}k_{2}^{b}e_{1})&=&r^{-6a-12b-3c-6d}s^{-6a-6b-3c-3d}
k_{1}^{a+2c}k_{2}^{b+2d}e_{2}^{2}e_{1};\\
(k_{1}^{c}k_{2}^{d}e_{2})(k_{1}^{a}k_{2}^{b}e_{1})(k_{1}^{c}k_{2}^{d}e_{2})&=&r^{-3a-6b-2c-3d}s^{-3a-3b-c}
k_{1}^{a+2c}k_{2}^{b+2d}e_{2}e_{1}e_{2};\\
(k_{1}^{a}k_{2}^{b}e_{1})^{2}(k_{1}^{c}k_{2}^{d}e_{2})^{2}&=&r^{-c+6d}s^{c+3d}
k_{1}^{a+2c}k_{2}^{b+2d}e_{1}e_{2}^{2}.
\end{eqnarray*}
\end{lem}
\qed

Now we are ready to prove one of the main results of this paper, which describes the group of algebra automorphisms of the algebra Hopf $\V$. Namely, we 
have the following
\begin{thm}
Let $\theta \in Aut_{\C}(\V)$ be an algebra automorphism of the Hopf algebra $\V$. Then for $l=1,2$, we have the following
\[
\theta(k_{l})=\lambda_{l}k_{l},\quad \theta(e_{1})=\gamma_{1}k_{1}^{a}K_{2}^{b}e_{1}, \quad \theta(e_{2})=\gamma_{2}k_{1}^{c}k_{2}^{d}e_{2}
\]
where $\lambda_{l}, \gamma_{l} \in \C^{\ast}$ and $a, b, c, d \in \Z$ such that $c=3b, a+3b+d=0$.
\end{thm}
{\bf Proof:} Let $\theta$ be an algebra automorphism of $\V$ and suppose that 
\[
\theta(e_{1})=\gamma_{1}k_{1}^{a}k_{2}^{b}e_{1},\quad
\theta(e_{2})=\gamma_{2}k_{1}^{c}k_{2}^{d}e_{2}.
\]

Note that the generators $e_{1}, e_{2}$ satisfy the following two-parameter quantum Serre relations
\begin{eqnarray*}
e_{2}^{2}e_{1}-(r^{-3}+s^{-3})e_{2}e_{1}e_{2}+r^{-32}s^{-3}e_{1}e_{2}^{2}=0;\\
e_{1}^{4}e_{2}-(r+s)(r^{2}+s^{2})e_{1}^{3}e_{2}e_{1}+rs(r^{2}+rs+s^{2})e_{1}^{2}e_{2}e_{1}^{2}\\-r^{3}s^{3}(r+s)(r^{2}+s^{2})e_{1}e_{2}e_{2}^{3}+r^{6}s^{6}e_{2}e_{1}^{4}=0.
\end{eqnarray*}

Since $\theta$ is an algebra automorphism of the Hopf algebra $\V$, we know that $\theta$ preserves the quantum Serre relations. In particular, we can derive the following system of 
equations via using the previous lemmas.
\begin{eqnarray*}
6a+18b+4c+12d&=&3a+12b+3c+9d;\\
6b+2c+6d&=&3a+12b+3c+9d;\\
-3a+c+3d&=&3a+12b+3c+9d;\\
-6a-6b&=&3a+12b+3c+9d;\\
12a+18b+8c+12d&=&9a+15b+6c+9d;\\
6a+12b+4c+6d&=&9a+15b+6c+9d;\\
3a+9b+2c+3d&=&9a+15b+6c+9d;\\
6b&=&9a+15b+6c+9d;\\
-6a-12b-3c-6d&=& -3a-6b-2c-3d;\\
-c&=&-3a-6b-2c-3d;\\
-6a-6b-3c-3d&=& -3a-3b-c;\\
c+3d&=&-3a-3b-c.\\
\end{eqnarray*}

Solving the previous system of equations, we shall obtain the following system of equations
\begin{eqnarray*}
3b&=&c;\\
a+c+d &=&0.
\end{eqnarray*}
Therefore, we have proved the theorem as desired.
\qed

\section{Hopf algebra automorphisms of $\V$}

In this section, we will determine all the Hopf algebra automorphisms of the 
Hopf algebra $\V$. We denote by $Aut_{Hopf}(\V)$ the group of all Hopf algebra automorphisms of $\V$. In particular, we shall prove that $Aut_{Hopf}(\V)$ is isomorphic 
to a torus of rank $2$.

To finish the task of this section, we need  to establish the following result 
\begin{thm}
Let $\theta \in Aut_{Hopf}(\V)$. Then for $l=1,2$, we have the following
\[
\theta(k_{l})=k_{l},\quad \theta(e_{l})=\gamma_{l}e_{l},
\]
for some $\gamma_{l}\in \C^{\ast}$. In particular, we have 
\[
Aut_{Hopf}(\V)\cong (\C^{\ast})^{2}.
\]
\end{thm}
{\bf Proof:} Let $\theta \in Aut_{Hopf}(\V)$ denote a Hopf algebra automorphism of $\V$, then we have $\theta \in Aut_{\C}(\V)$. According to {\bf Theorem 2.1.}, we shall have the following
\begin{eqnarray*}
\theta(k_{l})=\lambda_{l}k_{l};\\
\theta(e_{1})=\gamma_{1}k_{1}^{a}k_{2}^{b}e_{1};\\
\theta(E_{2})=\gamma_{2}k_{1}^{c}k_{2}^{d}e_{2};
\end{eqnarray*} 
for some $\lambda_{l}, \gamma_{l}\in \C^{\ast}$ for $l=1,2$, and $a,b,c,d \in \Z$ such that $3b=c, a+c+d=0$. 

First of all, we need to show that we have $\lambda_{l}=1$ for $l=1,2$. Since $\theta$ is a
Hopf algebra automorphism, we shall have the following
\[
(\theta \otimes \theta)(\Delta(k_{l}))=\Delta(\theta(k_{l}))
\]
for $l=1,2$, which implies the following
\[
\lambda_{l}^{2}=\lambda_{l}
\]
for $l=1,2$. Thus, we have $\lambda_{l}=1$ for $l=1,2$ as desired. 

Second of all, we need to show that we have $a=b=c=d=0$. Note that we have the following
\begin{eqnarray*}
\Delta(\theta(e_{1}))&=&\Delta(\gamma_{1}k_{1}^{a}k_{2}^{b}e_{1})\\
&=&\Delta(\gamma_{1}k_{1}^{a}k_{2}^{b})\Delta(e_{1})\\
&=&\gamma_{1}(k_{1}^{a}k_{2}^{b}\otimes k_{1}^{a}k_{2}^{b}) (e_{1}\otimes 1+k_{1}^{2}k_{2}^{-1}\otimes e_{1})\\
&=& \gamma_{1}k_{1}^{a}k_{2}^{b}e_{1}\otimes k_{1}^{a}k_{2}^{b}+\gamma_{1}k_{1}^{a}k_{2}^{b}k_{1}^{2}k_{2}^{-1}\otimes k_{1}^{a}k_{2}^{b}e_{1}\\
&=&\theta(e_{1})\otimes k_{1}^{a}k_{2}^{b}+k_{1}^{a}k_{2}^{b}k_{1}^{2}k_{2}^{-1}\otimes \theta(e_{1}).
\end{eqnarray*}

In addition, we also have the following
\begin{eqnarray*}
(\theta\otimes \theta)(\Delta(e_{1}))&=&(\theta\otimes \theta)(e_{1}\otimes 1+k_{1}^{2}k_{2}^{-1}\otimes e_{1})\\
&=& \theta(e_{1})\otimes 1+\theta(k_{1}^{2}k_{2}^{-1})\otimes \theta(e_{1})\\
&=& \theta(e_{1})\otimes 1+ k_{1}^{2}k_{2}^{-1}\otimes \theta(e_{1}).
\end{eqnarray*} 

Since we have $\Delta(\theta(e_{1}))=(\theta\otimes \theta)\Delta(e_{1})$, we shall have $a=b=0$. Note that we have $3b=c$ and $a+c+d=0$, thus 
we have $a=b=c=d=0$ as desired.

Conversely, it is obvious to see that the algebra automorphism $\theta$ defined by $\theta(k_{l})=k_{l}$ and $\theta(e_{l})=\gamma_{l}e_{l}$ for $l=1, 2$ is 
a Hopf algebra automorphism of $\V$. So we have proved the theorem.
\qed


\begin{thebibliography}{99999999}
\frenchspacing
\bibitem{BGH} Bergeron, N., Gao, Y. and Hu, N., Drinfeld doubles and 
Lusztig's symmetries of two-parameter quantum groups, {\it J. Algebra} {\bf
301} (2006), no. 1, 378--405.

\bibitem{BW1} Benkart, G. and Witherspoon, S., Two-parameter quantum groups and Drinfeld doubles, {\it Algebr. Represent. Theory,} {\bf 7 }(2004), 261--286.

\bibitem{BW2} Benkart, G. and Witherspoon, S., Representations of two-parameter quantum groups and Schur--Weyl duality, Hopf algebras, Lecture Notes in Pure and Appl. Math., {\bf 237}, pp. 65--92, Dekker, New York, 2004.

\bibitem{F} Fleury, O., Automorphisms of $\check{U}_{q}(\mathfrak b^{+})$, {\it Beitr\"{a}ge Algebra and Geom.}, {\bf Vol 38(2)} (1997), 343--356.

\bibitem{HQ} Hu, N.H.,  Shi, Q., The two-parameter quantum group of exceptional type $G_2$ and Lusztig's symmetries, {\it Pacific J. Math.},{\bf 230} (2007), no. 2, 327--345.

\bibitem{HW} Hu, N.H., Wang, X.L.,  Convex PBW-type Lyndon basis and restricted two-parameter quantum groups of type $G_2$,{\it Pacific J. Math.},{\bf 241} (2009), no. 2, 243--273.

\bibitem{T1} Tang, X., (Hopf) algebra automorphisms of the Hopf algebra ${\check U}^{\geq 0}_{r,s}({sl_{3}})$, submitted.
\bibitem{T2} Tang, X., Derivations of the two-parameter quantized enveloping algebra $U^{+}_{r,s}(B_{2})$, preprint.
\end{thebibliography}
\end{document}